\documentclass{amsart}

\usepackage{epsfig}
\usepackage{latexsym}
\usepackage{amsmath}
\usepackage{mathrsfs}
\newtheorem{theorem}{Theorem}[section]

\newtheorem{lemma}[theorem]{Lemma}

\newcommand{\cP}{\mathcal P}

\newcommand{\C}{{\mathscr C}}

\newcommand{\old}[1]{{}}

\title{A Hall-type theorem for triplet set systems based on medians in trees}

\author{Andreas Dress and Mike Steel}

\thanks{We thank the New Zealand Marsden Fund for supporting this work}

\address{AD: CAS-MPG Partner Institute for Computational Biology,\\
320 Yue Yang Road, 200031 Shanghai, China.\\
MS: Biomathematics Research Centre,\\
University of Canterbury, Christchurch, NZ\\
 }

\email{andreas@picb.ac.cn; m.steel@math.canterbury.ac.nz}



\begin{document}
\begin{abstract}
Given a collection $\C$ of subsets of a finite set $X$, let $\bigcup \C = \cup_{S \in \C}S$.
Philip Hall's celebrated theorem \cite{hall} concerning `systems of distinct representatives'  tells us that for any collection $\C$  of subsets of $X$ there exists an injective (i.e. one-to-one) function $f: \C \rightarrow X$ with
$f(S) \in S$ for all $S \in \C$ if and and only if $\C$ satisfies the property that for all non-empty subsets $\C'$ of $\C$ we have $|\bigcup \C'| \geq |\C'|$.  Here we show that if the condition
   $|\bigcup \C'| \geq |\C'|$ is replaced by  the stronger condition $|\bigcup \C'| \geq |\C'|+2$, then we obtain a characterization of
   this condition for a collection of 3-element subsets of $X$  in terms of
   the existence of an injective function from $\C$ to the vertices of a tree whose vertex set includes $X$ and that satisfies a certain median condition.
   We then describe an extension of this result to collections of arbitrary-cardinality subsets of $X$. 
\end{abstract}

\maketitle

\section{First result}

Given a tree $T=(V,E)$ and a subset $S$ of $V$ of size $3$, say $S = \{x,y,z\}$, consider the path in $T$ connecting $x,y$,  the path connecting $x,z$ and the path connecting $y,z$.  There is a unique vertex
that is shared by these three paths, the {\em median vertex} of $S$ in $T$, denoted ${\rm med}_T(S)$.
\begin{theorem}
\label{mainthm}
Let $X$ be a finite set, and suppose that $\C \subseteq \binom{X}{3}$, and $\bigcup \C = X.$ The following are equivalent:
\begin{itemize}
\item[(1)]
There exists a tree $T = (V,E)$ with $X \subseteq V$ for which the function $S \mapsto~{\rm med}_T(S)$ from $\C$ to $V$   is injective.
\item[(2)]
There exists a tree $T = (V,E)$  with  $X$ as its set of leaves, and all its other vertices of degree $3$,  for which the function $S \mapsto {\rm med}_T(S)$ from $\C$ to the set of interior vertices of $T$   is injective.
\item[(3)]
$\C$ satisfies the following property. For all non-empty subsets $\C'$ of $\C$ we have:
\begin{equation}
\label{ineq}
|\bigcup \C'| \geq |\C'|+2.
\end{equation}
\end{itemize}
\end{theorem}
In order to establish Theorem~\ref{mainthm} we first require a lemma.

Recall, from \cite{patch}, that a collection $\cP$ of subsets of a set $M$ forms a {\em patchwork} if it satisfies the following property:
$$A,B \in \cP \mbox{ and } A \cap B \neq \emptyset \Longrightarrow A \cap B, A \cup B \in \cP.$$

\begin{lemma}
\label{useful}
Let $X$ be a finite set, and suppose that $\C \subseteq \binom{X}{3}$, and $\bigcup \C = X.$  If $\C$ satisfies the condition described in Part (3) of Theorem~\ref{mainthm} then the collection $\cP$ of non-empty subsets
$\C'$ of $\C$ that satisfy $|\bigcup \C'| = |\C'|+2$ forms a patchwork.
\end{lemma}
{\em Proof:}
Suppose $\C_1, \C_2 \subseteq \C$, and that $\C_1 \cap \C_2 \neq \emptyset$.
Consider $$K: = |\bigcup (\C_1 \cap \C_2)| +  |\bigcup (\C_1 \cup \C_2)| .$$ By (\ref{ineq}) we have:
\begin{equation}
\label{patching1}
K \geq (|\C_1 \cap \C_2|+2)+(|\C_1 \cup \C_2|+2) = |\C_1|+|\C_2|+4,
\end{equation}
and we also have:
\begin{equation}
\label{patching2}
K  \leq   |(\bigcup \C_1) \cap (\bigcup \C_2)| + |(\bigcup \C_1) \cup (\bigcup \C_2)|        = |\bigcup \C_1|+|\bigcup \C_2|.
\end{equation}

Notice that the right-hand term in (\ref{patching1}) and in (\ref{patching2}) are equal,  since
$|\bigcup \C_i| = |\C_i|+2$ as $\C_i \in \cP$ for $i=1,2$, and
thus, the inequality in ({\ref{patching1}) is an equality, and so $|\bigcup (\C_1 \cap \C_2)|  = |\C_1 \cap \C_2|+2$ and  $|\bigcup (\C_1 \cup \C_2)|  = |\C_1 \cup \C_2|+2$, as required.
\hfill$\Box$

{\em Proof of Theorem~\ref{mainthm}:}
The implication (2) $\Rightarrow$ (1) is trivial. For the reverse implication suppose that $T$ satisfies the property described in (2).  First delete from $T$ any vertices and edges that are not on a path between two vertices in $X$.  Next attach to every interior  (non-leaf) vertex $v \in X$ a new edge for which the adjacent new leaf is assigned the label $x$, and henceforth do not regard  $v$ as an element of $X$. Next replace each maximal path of degree $2$ vertices
by a single edge. Finally replace each vertex $v$ of degree $d>3$ by an arbitrary tree that has $d$ leaves that we identify with the neighboring vertices of $v$ and whose remaining vertices have degree 3. These four processes result in a tree $T'$ that has $X$ as its set of leaves, and which has all its remaining vertices of degree $3$ (i.e. a `binary phylogenetic $X$--tree' \cite{sem})  and for which
the median vertices of the elements of $\C$ remain distinct. Thus (1) and (2) are equivalent.

Next we show that (2) $ \Rightarrow$ (3).  Suppose $T$ satisfies the condition (2) and that $\C'$ is a non-empty subset of $\C$.  Consider the minimal subtree of $T$ that connects the leaves in
$\bigcup \C'$.  This tree has at least $|\C'|$ interior vertices that are of degree $3$. However by a simple counting argument, any tree that has $k$ interior vertices of degree  $3$ must have at least $k+2$ leaves, and so (\ref{ineq}) holds.

The remainder of the proof is devoted to establishing that (3) $\Rightarrow$ (2).  We use induction on $n:=|X|$.  The result clearly holds for $n=3$, so suppose it holds whenever $|X| < n, n \geq 4$ and that $X$ is a set of size $n$.  For $x \in X$, let $n_\C(x)$ be the number of triples in $\C$ that contain $x$.  If there exists $x \in X$ with $n_\C(x) =1$, then select the unique triple in $\C$ containing $x$, say
$\{a,b,x\}$ and let $X' = X-\{x\}, \C' = \C-\{\{a,b,x\}\}.$  Then $\bigcup \C' = X'$ and  $\C'$ satisfies condition (\ref{ineq}) and so, by induction, there is a tree $T'$ with leaf set $X'$ for which the median vertices of elements in $\C'$
are all distinct vertices of $T'$.  Let $T$ be the tree obtained from $T'$ by subdividing one of the edges in the path in $T'$ connecting $a$ and $b$ and making the newly-created vertex of degree 2 adjacent to $x$ by a new edge. Then $T$ satisfies the requirements of Theorem~\ref{mainthm}(2), and thereby establishes the induction step in this case.

Thus we may suppose that $n_\C(x)>1$ holds for all $x \in X$. In this case we claim that there exists
$x \in X$ with $n_\C(x)=2$.   Let us count the set $\Omega:= \{(x, S): x \in S \in \C\}$ in two different ways.
We have
\begin{equation}
\label{eq1}
|\Omega| = \sum_{x \in X} n_\C(x) \geq 2k+3(n-k),
\end{equation}
where $k = |\{x\in X: n_\C(x) =2\}|.$
On the other hand,
\begin{equation}
\label{eq2}
|\Omega| = 3|\C| \leq 3(n-2),
\end{equation}
where the latter inequality follows from Inequality  (\ref{ineq}) applied to $\C' =\C$.
Combining (\ref{eq1}) and (\ref{eq2}) gives $2k+3(n-k) \leq 3n-6$, and so $k \geq 6$. Thus, since $k>0$, there exists $x \in X$ with $n_\C(x) =2$, as claimed.

For any such $x \in X$ with $n_\C(x)=2$, let $\{a,b,x\}, \{a',b', x\}$ be the two elements of $\C$ containing $x$.  Without loss of generality there are two cases:
\begin{itemize}
\item[(i)]  $a=a', b \neq b'$; or
\item[(ii)] $\{a,b\} \cap \{a',b'\} = \emptyset$.
\end{itemize}

In case (i) let $$X' := X-\{x\}, \C':= \C-\big\{\{a,b,x\}, \{a, b',x\}\big\}, \C_1 := \C' \cup \big\{\{a,b,b'\}\big\}.$$
Note that $\bigcup \C_1 = X'$.  Suppose that $\C_1$ fails to satisfy the condition described in Part (3) of Theorem~\ref{mainthm}.  Then there is a
subset of $\C_1$ that violates Inequality (\ref{ineq})  of the form $\C^1 \cup \big\{\{a,b,b'\}\big\}$ where $a,b,b' \in \bigcup \C^1$ and
$\C^1 \subseteq \C'$.  But in that case  $\C^1 \cup \big\{\{a,b,x\}, \{a,b',x\}\big\}$ would violate Inequality (\ref{ineq}), which is impossible since
Inequality  (\ref{ineq}) applies to this set, being a non-empty subset of $\C$.  Thus, $\C_1$ satisfies Part (3) of Theorem~\ref{mainthm}.
Since $\bigcup \C_1 = X'$, which has one less element than $X$, the inductive hypothesis furnishes a tree
$T'$ with leaf set $X'$ that satisfies the requirements of Theorem~\ref{mainthm}(2).  Now consider the edge of $T'$ that is incident with leaf $b'$. Subdivide this edge and make the newly created midpoint vertex adjacent to a leaf labelled $x$. This gives a tree $T$ that has $X$ as its set of leaves, and with all its interior vertices of degree 3; moreover  the medians of the elements
of $\C$ are all distinct (note that the median of $\{x,a,b'\}$ is the newly-created vertex adjacent to $x$, while the median of $\{x,a,b\}$  corresponds to the median vertex of $\{a,b,b'\}$  in $T'$ and therefore is a different vertex in $T$ to any other median vertex of an element of $\C$).

In case (ii), let $$X' := X-\{x\}, \C':= \C-\big\{\{a,b,x\}, \{a',b',x\}\big\},$$
and $$ \C_1:= \C' \cup \{\{a,a',b\}\}, \C_2 :=  \C' \cup \big\{\{a,a',b'\}\big\}.$$
Note that $\bigcup \C_1 = \bigcup \C_2 = X'$.  We will establish the following

{\bf Claim:} One or both of $\C_1$ or $\C_2$ satisfies the condition described in Part (3) of Theorem~\ref{mainthm}.

\noindent Suppose to the contrary that
both sets fail  the condition described in Theorem~\ref{mainthm}(3). Then there is a  subset of $\C_1$ that violates Inequality (\ref{ineq}),  and it must be of the form
$\C^1 \cup  \{\{a,a',b\}\}$ where $\C^1 \subseteq \C'$, $a,a',b \in \bigcup \C^1$ and $b' \not\in \bigcup \C^1$ (the last claim is justified by the observation that if $b' \in \bigcup \C^1$ then $\C^1 \cup \big\{\{a,b,x\}, \{a',b',x\}\big\}$ would violate the condition
described in Part (3) of Theorem~\ref{mainthm}).
Similarly  a subset of $\C_2$ that violates Inequality (\ref{ineq})  is of the form
$\C^2 \cup  \big\{\{a,a',b'\}\big\}$ where $\C^2 \subseteq \C'$, $a,a',b' \in \bigcup \C^2$ and $b \not\in \bigcup \C^2.$
Now, let $\cP$ be the subset of $\C$ defined in the statement of  Lemma~\ref{useful}.  Then the sets
$$\C_1: = \C^1 \cup \big\{\{x,a,b\}, \{x,a',b'\} \big\};  \mbox{ and }  \C_2: = \C^2 \cup \big\{\{x,a,b\}, \{x,a',b'\} \big\}$$
are both elements of $\cP$ and they have non-empty intersection, since they both contain $\{x,a,b\}$ (indeed they also share $\{x,a',b'\}$). Thus,  Lemma~\ref{useful} ensures that $\C_1 \cap \C_2$ is also an element of
$\cP$.  However $\C_1 \cap \C_2$ is of the form $\C^3 \cup \big\{\{x,a,b\}, \{x,a',b'\} \big\}$ where $\C^3 \subseteq \C'$, and neither $x$, nor $b$, nor $b'$ are elements of $\bigcup \C^3$ because by our choice of $x$, $x$ only occurs in the two triples $\{x,a,b\}$ and $\{x, a', b'\}$, and because $b'  \not\in \bigcup \C^1$ and
$b \not\in \bigcup \C^2.$ Since $\C^3$ is a subset of $\C$, $\C^3$ satisfies Inequality  (\ref{ineq}) which implies that (\ref{ineq}) must be a strict inequality for $\C_1 \cap \C_2$, contradicting our assertion that
$\C_1 \cap \C_2 \in \cP$. This justifies our claim that  either $\C_1$ or $\C_2$ satisfies the part (3) of Theorem~\ref{mainthm}.

We may suppose then, without loss of generality, that $\C_1$ satisfies part (3) of Theorem ~\ref{mainthm}. Since $\bigcup \C_1 = X'$, which has one less element than $X$, the inductive hypothesis furnishes a tree
$T'$ with leaf set $X'$ that satisfies the requirements of Theorem~\ref{mainthm}(2).  Now,  consider the edge of $T'$ that is incident with leaf $a'$. Subdivide this edge and make the newly created midpoint vertex adjacent to a leaf labelled $x$ by a new edge.  This gives a tree $T$ that has $X$ as its set of leaves, and with all vertices of degree 3; moreover, regardless of where $b'$ attaches in $T$ the medians of the elements
of $\C$ are all distinct (note that the median of $\{x,a',b'\}$ is the newly-created vertex adjacent to $x$, while the median of $\{x,a,b\}$ corresponds to  the median vertex of  $\{a,a',b\}$ in $T'$  and therefore is a different vertex in $T$ to any
other median vertex of an element of $\C$).   This completes the proof.
\hfill$\Box$

\section{An extension}

For a subset $Y$ of $X$ of size at least $3$, and a tree $T = (V,E)$, with $X \subseteq V$, let $${\rm med}_T(Y) := \{{\rm med}_T(S): S \subseteq Y, |S|=3\}.$$  Thus, ${\rm med}_T(Y)$ is a subset of the vertices of $T$, moreover if $X$ is the set of leaves of $T$, then ${\rm med}_T(Y)$ is a subset of the interior vertices of $T$.
\begin{theorem}
\label{mainthm2}
Let $X$ be a finite set, and suppose that $\C$ is a collection of subsets of $X$, each of size at least $3$, and with
 $\bigcup \C = X$. 
 The following are equivalent:
\begin{itemize}
\item[(1)]
There exists a tree $T = (V,E)$  with  $X$ as its set of leaves, and all its other vertices of degree $3$,  for which $\{{\rm med}_T(Y): Y \in \C\}$ is a partition of the set of interior vertices of $T$.
\item[(2)]
$\C$ satisfies the following property. For all non-empty subsets $\C'$ of $\C$ we have:
\begin{equation}
\label{inex}
|\bigcup \C'|-2 \geq \sum_{Y \in \C'}(|Y|-2),
\end{equation}
and this last inequality is an equality when $\C'=\C$. 
\end{itemize}
\end{theorem}
{\em Proof:} We first show that (1) $\Rightarrow$ (2). Select a tree $T$ satisfying the requirements of  Part (1) of Theorem~\ref{mainthm2}. For a nonempty  subset $\C'$ of $\C$, the minimal subtree $T'$ of $T$ connecting the leaves in $\bigcup \C'$ has $k:=|\bigcup \C'|$ leaves, and $k-2$ vertices that are of degree $3$.
By the partitioning assumption, each element $Y \in \C'$ generates $|Y|-2$ median vertices in $T$ and these sets of median vertices are pairwise disjoint for different choices of $Y \in \C'$. Moreover, distinct interior vertices of $T$  correspond to different degree 3 vertices in $T'$, and so the number of degree $3$ vertices in $T'$ can be no smaller than the sum of $|Y|-2$ over all $Y \in \C'$. This establishes Inequality~(\ref{inex}). For the case where $\C'= \C$, note that $T$ has $\bigcup \C = X$ as its leaf set, and by the partitioning assumption, each of its $|X|-2$ interior vertices occurs in one set ${\rm med}_T(Y)$ for some $Y \in \C$ and so $|X|-2 \leq \sum_{Y \in \C'}(|Y|-2)$ which, combined with (\ref{inex}), provides the desired equality.

To show (2) $\Rightarrow$ (1) select for each set $Y \in \C$ a collection $\C_Y$ of 3-element subsets of $X$ of cardinality $|Y|-2$  for which $\bigcup C_Y = Y$, and that satisfies that condition that
for every nonempty subset $\C'$ of $\C_Y$  we have 
$\bigcup \C' \geq |\C'|+2$; such a selection is straightforward - for example, if 
$Y = \{y_1,\ldots, y_m \}$ then we can take
\begin{equation}
\label{spec}
\C_Y = \big\{\{y_1,y_2,y_3\}, \{y_1, y_2, y_4\}, \ldots, \{y_1,y_2, y_m \} \big\}.
\end{equation}
 We first establish the following:

\noindent{\bf Claim:} $\C_*: = \cup_{Y \in \C} \C_Y$  is a collection of 3-element subsets of $X$ that satisfies Inequality~(\ref{ineq}) in Theorem~\ref{mainthm}.

To see this, suppose to the contrary that there exists a subset $\C''$ of $\C_*$ for which Inequality (\ref{ineq}) fails. Write $\C'' = S_1 \cup S_2 \cup \cdots \cup S_k$ where $1 \leq k \leq |\C|$, and where $S_i$  is a nonempty set of $3$--element subsets of $X$ that are selected from the same set  (let us call it $Y_i$)  from $\C$  (note that the fact that $|Y_1\cup Y_2| - 2 \ge |Y_1| - 2 + |Y_2| - 2$ must hold for all $Y_1, Y_2$ in $\C$ implies that $|Y_1| + |Y_2| - |Y_1\cap Y_2| - 2 \ge |Y_1| - 2 + |Y_2| - 2$ and, hence, $2\ge |Y_1\cap Y_2| $ must hold for all $Y_1, Y_2$ in $\C$).  By our assumption regarding the set of triples $\C''$ we have $|\bigcup \C''| \leq |\C''|+1$ and so, if we let $W_i := \bigcup S_i$ we have  $\bigcup \C'' = \bigcup_{i=1}^k W_i$, and consequently
\begin{equation}
\label{in1}
|\bigcup_{i=1}^k W_i|   \leq \sum_{i=1}^k |S_i| +1.
\end{equation}
For  $\C': = \{Y_1, \ldots, Y_k\} \subseteq \C$ we have
\begin{equation}
\label{in2}
|\bigcup \C'| \geq \sum_{i=1}^k(|Y_i|-2)+2 =  \sum_{i=1}^k |Y_i|-2k+2.
\end{equation}
On the other hand,
$$|\bigcup \C' |  \leq |\bigcup_{i=1}^k W_i| + \sum_{i=1}^k(|Y_i - W_i|)  = |\bigcup_{i=1}^k W_i| + \sum_{i=1}^k(|Y_i| - |W_i|),$$ 
since $W_i \subseteq Y_i$.  By the condition imposed on the construction of $\C_Y$ we have  $|W_i| \geq |S_i|+2$ for each $i$, and so, substituting this, and (\ref{in1})  into the previous inequality gives:
$$|\bigcup \C'| \leq  \sum_{i=1}^k |S_i| +1 + \sum_{i=1}^k |Y_i| - \sum_{i=1}^k(|S_i|+2) =\sum_{i=1}^k |Y_i|-2k+1, $$ 
which compared with (\ref{in2}) gives $1\geq 2$, a contradiction. This establishes that $C_*$ satisfies Inequality~(\ref{ineq}) in Theorem~\ref{mainthm}.

By Theorem~\ref{mainthm} it now follows that there is a tree $T=(V,E)$ with leaf set $X$ for which the function $S \mapsto {\rm med}_T(S)$ is injective from $\C_*$ to the set of interior vertices of $T$.
Now,  for $Y \in \C$, we have
\begin{equation}
\label{eqx}
{\rm med}_T(Y) =  \{{\rm med}_T(S): S \subseteq Y, |S|=3\} =   \{{\rm med}_T(S): S \in C_Y\}.
\end{equation}
The second equality in (\ref{eqx}) requires some justification.  Recalling our particular choice of $C_Y$ from (\ref{spec}), and noting that the medians of the triples in $C_Y$ are distinct vertices of $T$ it follows that $T|Y$ has the structure of a path connecting $y_1, y_2$ with each of the remaining  leaves $y \in Y-\{y_1, y_2\}$ separated from this path by just one edge.  Consequently if a vertex $v$ of $T$ is the median of three
leaves in $Y$, then it is also the median of a triple $\{y_1, y_2, y\}$ for some $y \in Y-\{y_1, y_2\}$; that is, it is an element of $\{{\rm med}_T(S): S \in C_Y\}.$

Consequently,  $\{{\rm med}_T(Y): Y \in \C\}$  are disjoint subsets of the set of interior vertices of $T$.
Moreover each interior vertex of $T$ is covered by  $\{{\rm med}_T(Y): Y \in \C\}$  since the number of interior vertices is $|X|-2$ and by assumption $|X|-2 = \sum_{Y \in \C}(|Y|-2) = |\C_*|$.   This establishes the implication (2) $\Rightarrow$ (1) and thereby completes the proof.
\hfill$\Box$

\end{document}